\documentclass[amsmath,article,aps,fleqn]{revtex4}

\begin{document}

\title{A determinant formula for the partition function $p(7k+a)$ }
\author{Jerome Malenfant}
\affiliation{American Physical Society\\ Ridge, NY}
\date{\today}
\begin{abstract} We derive expressions for the partition function $p(n)$, with $n$ in the form $7k+a$,
as $(k+1)$-dimensional determinants.   
\end{abstract}

\maketitle
\section{Introduction and Main Result}
In a recent article \cite{Mal}, a formula was derived for the partition function $p(n)$ as the $ (n+1)$-dimensional determinant   
           \begin{eqnarray}
          p(n) &=&  \left| \begin{array}{rrrrrrrrrr}
                       1& ~ & ~&~&~&~&~&~&~&1\\
            -1& 1& ~ & ~&~&~&~&~&~&0\\
            -1& -1& ~1& ~ & ~&~&~&~&~&0\\
            0& -1& -1& ~1 & ~ & ~&~&~&~&0\\
            0 & 0 & -1& -1&~1 & ~ &~&~& ~&0\\
            1& 0 & 0& -1&   -1&~1& ~ &~&~& 0\\
            0&  1& 0 & 0& -1& -1&~1  &  ~ &~&0\\
             1&0&  1& 0 & 0& -1& -1&~ &  ~ &0\\   
            ~ \vdots ~& ~&~&~&~&~&~&~~\ddots &~& \vdots 
            \end{array} \right|_{(n+1) \times (n+1)},
           \end{eqnarray}            
using the generating function
     \begin{eqnarray}
     \sum_{n=0}^{\infty} p(n) q^n = \frac{1}{(q)_{\infty} }, ~~ {\rm where}~~(q)_{\infty} \equiv  \prod_{k=1}^{\infty} (1-q^k).
     \end{eqnarray}
The coefficients along the diagonals  of the above matrix  are from the expansion of $(q)_{\infty}$ in powers of $q$, (sequence A010815 in OEIS \cite{integer}).  In Appendix A of that article, $(k+1)$-dimensional determinant formulas for the partition functions $p(5k+4), ~p(7k+5),$ and $p(25k+24)$ were  found  by using some identities due to  Ramanujan \cite{Ram}.  The $p(5k+4)$ and $p(25k+24)$ expressions were then extended to $p(5k+a_5)$ and $p(25k+a_{25})$
for $a_5 =0,1,2,3$ and $a_{25} = 4,9,14,19$, respectively.

In this article we will continue this development and generalize the $p(7k+5)$ result,
    \begin{eqnarray}
          p(7k+5) &=&  \left| \begin{array}{rrrrrrrrrr}
             1& ~ & ~&~&~&~&~&~&~&~7\\
            -8& 1& ~ & ~&~&~&~&~&~&~21\\
            20& -8& 1& ~ & ~&~&~&~&~&~14\\
            ~0& ~20& -8& 1 & ~ & ~&~&~&~&~56\\
            -70 & ~0 & ~20& -8&1 & ~ &~&~& ~&-35\\
            ~64& -70 & ~0& ~20&   -8&1& ~ &~&~& -28\\
            ~56&  ~64& -70 & ~0& ~20& -8&~1  &  ~ &~&-70\\
             ~ 0&~56&  ~64& -70 & ~0& ~20& -8&~ &  ~ &~35\\   
             ~\vdots ~& ~&~&~&~&~&~&~~\ddots &~&  \vdots ~
            \end{array} \right|_{(k+1) \times (k+1)} ,
           \end{eqnarray}
which was derived from            
           \begin{eqnarray}
          \sum_{k=0}^{\infty} p(7k+5)q^k &=& \frac{(q^7)^7_{\infty}}{(q)^8_{\infty} }
           \left[ 7~ \frac{(q)^4_{\infty} }  {(q^7)^4_{\infty}}  +49q \right]  
           \end{eqnarray}  
using the expansions
          \begin{eqnarray}
          (q)^8_{\infty} &=& 1 -8q +20q^2 -70q^4 + 64 q^5 + \cdots,  \nonumber \\
          7 (q^7)^3_{\infty}(q)^4_{\infty} + 49q (q^7)^7_{\infty} &=& 7 +21q +14q^2 +56q^3 -35q^4 - \cdots,  \nonumber
           \end{eqnarray}     
to $p(7k+a_7)$ for $a_7 = 0, \ldots, 6$.  
We will follow Ramanujan \cite{Ram} and make the replacement $q \rightarrow q^{1/7}$ in eq. (2). $(q^{1/7})_{\infty}$  expands as the sum $\sum b_n q^{n/7}$, for which the only nonzero $ b_n$ are for $n \equiv 0, 1, 2, {\rm or~} 5\mod 7 $.   Furthermore, the  $n \equiv  2 \mod 7 $ terms sum to $-(q^7)_{\infty} q^{2/7}$.  We can  therefore write,      
           \begin{eqnarray}
          \frac{ (q^{1/7})_{\infty}} {(q^7)_{\infty}} = J_1+q^{1/7}J_2 -q^{2/7} +q^{5/7}J_3.
          \end{eqnarray}
The $J$'s are power series in $q$ with integer exponents and satisfy the identities \cite{Ram}
\footnote{Note that in ref. 3 there is a misprint in the exponent of the last term in eq. (24.4), corresponding to eq. (6c);  $57q^3$ should be replaced in (24.4) by $57q^2$.}
           \begin{subequations}
          \begin{eqnarray}
            J_1J_2J_3 &=& -1 ;\\  \frac{J_1^2}{J_3} + \frac{J_2}{J_3^2}&=& q;  \\
           J_1^7 +J_2^7q + J_3^7 q^5 &=& \frac{(q)^8_{\infty}}{(q^7)_{\infty}^8} + 14q \frac{(q)^4_{\infty}}{(q^7)_{\infty}^4}
           +57q^2;\\
           J_1^3J_2 +J_2^3J_3q + J_1J_3^3 q^2 &=& - \frac{(q)^4_{\infty}}{(q^7)_{\infty}^4}-8q ;\\ 
            J_1^2J_2^3 +J_1^3J_3^2q +J_2^2J_3^3 q^2  &=& - \frac{(q)^4_{\infty}}{(q^7)_{\infty}^4} -5q.  
           \end{eqnarray}
           \end{subequations}
We have  then
          \begin{eqnarray}
          \sum_{n=0}^{\infty} p(n) q^{n/7} = \frac{1} {(q^7)_{\infty}}~ \frac{1} { J_1+q^{1/7}J_2 -q^{2/7} +q^{5/7}J_3}.
          \end{eqnarray}
As per Ramanujan, we now multiply and divide $ (J_1+q^{1/7}J_2 -q^{2/7} +q^{5/7}J_3)$ by the product
            \begin{eqnarray*}     
            \prod_{n=1}^6  \left(J_1+\omega^n q^{1/7}J_2 -\omega^{2n}q^{2/7} +\omega^{5n} q^{5/7}J_3\right)  
             \end{eqnarray*}        
where $ \omega \equiv e^{2 \pi i /7} $.  To simplify the notation and to avoid writing out a lot of subscripts, we 
define $   x \equiv  qJ_2^7/J_1^7 $ and $ a\equiv J_1/J_2^2$.  Then
          \begin{eqnarray}
        J_1+q^{1/7}J_2 -q^{2/7} +q^{5/7}J_3 =J_1 \left( 1 + x^{1/7} -x^{2/7}a -  x^{5/7}a^3 \right).
          \end{eqnarray}
Eq. (7) then becomes          
     \begin{eqnarray}
          \sum_{n=0}^{\infty} p(n) q^{n/7}  
          =  \frac{1} {(q^7)_{\infty}}~ \frac{ J_1^6  \prod_{n=1}^6  \left(1 + \omega^n x^{1/7}- \omega^{2n} x^{2/7} a
          -  \omega^{5n} x^{5/7}a^3 \right) } 
          {J_1^7  \prod_{n=0}^6  \left(1 + \omega^n x^{1/7}- \omega^{2n} x^{2/7}a -  \omega^{5n}  x^{5/7}a^3 \right) }
          \equiv  \frac{1} {(q^7)_{\infty}}~ \frac{N} { D}.
           \end{eqnarray}
The numerator $N$ and the denominator $D$ are expressible as  the sums
              \begin{eqnarray}
              N=  J_1^6\sum_{k=0}^{30} c_kx^{k/7} , ~~D = J_1^7  \sum_{k=0}^{5} d_kx^{k}  , 
              \end{eqnarray}
with coefficients $c_k$ and $d_k$, respectively.   Expanding the product in the denominator in eq. (9)  we get   
           \begin{subequations}                         
          \begin{eqnarray}     
        D  = J_1^7\left [~1 +x \left(1+7a+14a^2 -7a^4\right) -x^2 \left(8a^7 -14a^8\right) -14x^3 a^{11} -7x^4 a^{16} 
             -x^5 a^{21}~ \right] .
             \end{eqnarray}
Converting to $(q,J)$ notation, this is
             \begin{eqnarray}
              D = J_1^7 +q\left(J_2^7+7 J_1J_2^5+14J_1^2J_2^3 +7J_1^{5}J_3 \right)   
              - q^2 \left(8 -14J_1^3J_3^2 \right)  + 14 q^3 J_2^2 J_3^3  +7 q^4 J_2J_3^5  + q^5 J_3^7  ,
             \end{eqnarray}                   
which, after some manipulation and using the identities (6a-e) above, becomes 
             \begin{eqnarray}
             D = \frac{(q)^8_{\infty}}{(q^7)_{\infty}^8} .
             \end{eqnarray}
              \end{subequations} 
And so           
            \begin{eqnarray}
          \sum_{n=0}^{\infty} p(n) q^{n/7}  &=&  \frac{(q^7)_{\infty} ^7}{ (q)_{\infty}^8 }~ J_1^6 \sum_{k=0}^{30} c_k x^{k/7}
          \nonumber \\
          &=& \frac{(q^7)_{\infty} ^7}{ (q)_{\infty}^8 } \left( H_1+ H_2 q^{1/7} +  H_3q^{2/7}+  H_4 q^{3/7}+ 
                     H_5 q^{4/7} +  H_6 q^{5/7}   + H_7 q^{6/7}\right),
             \end{eqnarray}
where
             \begin{eqnarray*}
                     H_1 &=& ~~J_1^6 ~~ (~c_0+c_7x + c_{14}x^2 +c_{21}x^3 + c_{28}x^4~) ;\\
                     H_2 &=& J_1^5J_2 ~(~c_1 +c_8 x+c_{15}x^2 +c_{22}x^3 +c_{29}x^4~) ; \\
                     H_3 &=& J_1^4J_2^2~ ( ~c_2+ c_9x +c_{16}x^2 +c_{23}x^3 +c_{30}x^4~) ; \\
                     H_4 &=& J_1^3J_2^3~ (~ c_3 + c_{10}x +c_{17}x^2 +c_{24}x^3~);  \\
                     H_5 &=& J_1^2J_2^4~ (~ c_4+c_{11}x +c_{18}x^2 +c_{25}x^3 ~); \\
                     H_6 &=& J_1J_2^5~ (~ c_5+ c_{12}x +c_{19}x^2 +c_{26}x^3~ ) ; \\
                     H_7 &=& ~~J_2^6~~ ( ~c_6 +c_{13}x +c_{20}x^2 +c_{27}x^3~ ).                     
                 \end{eqnarray*}             
The generating-function equations are then 
          \begin{eqnarray}
          \sum_{k=0}^{\infty} p(7k+a_7) q^k =  \frac{(q^7)_{\infty} ^7}{ (q)_{\infty} ^8}H_{1+a_7}(q) ; ~~~a_7 = 0,1, \ldots ,6,
         \end{eqnarray}
 and we have
           \begin{eqnarray}
          p(7k+a_7) &=&  \left| \begin{array}{cccccccccc}
             ~1& ~ & ~&~&~&~&~&~&~& Z_0\\
            -8& ~1& ~ & ~&~&~&~&~&~& Z_1\\
            ~20& -8& ~1& ~ & ~&~&~&~&~& Z_2\\
            ~0& ~20& -8& ~1 & ~ & ~&~&~&~& Z_3\\
            -70 & ~0 & ~20& -8&~1 & ~ &~&~& ~& Z_4\\
            ~64& -70 & ~0& ~20&   -8&~1& ~ &~&~&  Z_5\\
            ~56&  ~64& -70 & ~0& ~20& -8&~1  &  ~ &~& Z_6\\
             ~ 0&~56&  ~64& -70 & ~0& ~20& -8&~ &  ~ & Z_7\\   
             ~\vdots & ~&~&~&~&~&~&~~\ddots &~& ~ \vdots 
            \end{array} \right|_{(k+1) \times (k+1)}
           \end{eqnarray}  
where the coefficients in the last column are from the expansions
            \begin{eqnarray}
            (q^7)_{\infty} ^7H_{1+a_7}(q) = Z_0 + Z_1q + Z_2q^2 + \cdots .
            \end{eqnarray}         
Eq. (3) gives the first few of these coefficients for $a_7 = 5$.
             
The $c_k$ coefficients can be found from the recurrence relation       
           \begin{eqnarray}
      c_n +c_{n-1} -ac_{n-2} -a^3c_{n-5} = \left \{ \begin{array} {l} d_{n/7} ~~{\rm if} ~ n \equiv 0\mod 7, \\ ~~0 ~~~~{\rm otherwise,}
          \end{array} \right.  
      \end{eqnarray}       
where the $d_k$  coefficients are given by eq. (11a).  The coefficient $c_0$ equals 1; the rest are listed below:

      \begin{eqnarray*} 
          \begin{array} {l} c_1\\c_2\\c_3 \\ c_4 \\ c_5 \\c_6\\ c_7\\ c_8 \\c_9 \\c_{10} \end{array}~~
          \begin{array}{r} -1 \\ 1+a \\ -1-2a \\ 1+3a+a^2 \\ -1-4a-3a^2+a^3 \\1+5a +6a^2 -a^3 \\
                 a+4a^2 -a^3 -5a^4\\ a^2 +6a^3 +2a^4\\  -a^3 -4a^5   \\ a^3 +2a^4 +3a^5 +a^6 \\ \end{array}
                 ~~~~~~~~\begin{array} {l} c_{11} \\c_{12} \\c_{13}\\c_{14} \\ c_{15} \\ c_{16} \\c_{17}\\ c_{18}\\ c_{19} 
                 \\c_{20} \end{array}~~  
                 \begin{array} {r}  2a^4 +3a^5 -6a^6 \\ 3a^5+8a^6 -4a^7 \\ a^6 \\ a^6 +6a^8 \\ 3a^7 -3a^8 +a^9  \\ 
                 6a^8 -a^9 \\ 6a^9 -3a^{10} \\ a^9 +2a^{10} \\4a^{10} +3a^{11} \\ -4a^{11} +a^{12}  \end{array}
                 ~~~~~~~~~\begin{array} {l} c_{21} \\c_{22} \\c_{23}\\c_{24} \\ c_{25} \\ c_{26} \\c_{27}\\ c_{28}\\ c_{29} 
                 \\c_{30} \end{array}~~ 
                 \begin{array}{r} a^{12} \\ a^{12} -2a^{13} \\ 5a^{13} \\ a^{14} \\ a^{15} \\ a^{15} \\ -a^{16} \\ 0 \\ 0 \\ a^{18} 
                 \end{array}             
           \end{eqnarray*}  
With these expressions, and using the identity (6b), the $H$ functions become
     \begin{eqnarray} 
     \nonumber\\
    H_1 &=& 2J_1^2J_2^8  +2J_1^3J_2^6 -J_1^4J_2^4 -13J_1^5J_2^2+ 11J_1^6; \nonumber  \\ \nonumber ~\\
    H_2 &=& 5J_1^2J_2^7 -9J_1^3 J_2^5 +15J_1^4J_2^3-15J_1^5J_2 -3J_1^7J_3; \nonumber\\~ \nonumber\\
    H_3    &=& 11 J_1^2J_2^6 -31J_1^3J_2^4 +26J_1^4J_2^2 -5J_1^5 +J_1^8J_3^2 ; \nonumber\\~ \nonumber\\
    H_4  &=& J_1J_2^7 +8J_1^2J_2^5 -18J_1^3J_2^3 +11J_1^4J_2 +5J_1^6J_3; \\~ \nonumber\\
    H_5 &=& 3J_1J_2^6 +3J_1^2J_2^4 -12J_1^3J_2^2 +12J_1^4 -J_1^7J_3^2 ; \nonumber\\~ \nonumber\\
    H_6  &=& 7 J_1J_2^5 -7J_1^2J_2^3 -14J_1^3J_2 -7J_1^5J_3; \nonumber\\ ~ \nonumber\\
    H_7  &=& J_2^6 + J_1J_2^4 +17J_1^2J_2^2 -10J_1^3 + 2J_1^6J_3^2. \nonumber\\ \nonumber
        \end{eqnarray}
$H_6$ simplifies  to the expression inside the brackets in eq. (4) using the additional identities (6c-e). 

 It remains to determine the expansions for the $J$ functions.  From eq. (5) we get     
    \begin{eqnarray}
    J_1(q) &=& \frac{1}{(q^7)_{\infty}}\left\{-1 + \sum_{k=0}^{\infty} q^{k(42k-1)} \left[ ~1 + q^{2k} +q^{14k+1} 
     - q^{30k+5} - q^{42k+10} - q^{44k+11}  - q^{56k+18}+q^{72k+30} ~ \right] \right\}, \nonumber  \\
   J_2(q) &=&  \frac{-1}{(q^7)_{\infty}}   \sum_{k=0}^{\infty} q^{k(42k+5)} \left[ ~1 + q^{14k+2} -q^{18k+3} 
    - q^{32k+8} - q^{42k+13}-q^{56k+22}  + q^{60k+25}  + q^{74k+37} ~ \right] , \\
   J_3(q) &=& \frac{1}{(q^7)_{\infty}}  \sum_{k=0}^{\infty} q^{k(42k+11)} \left[ ~1 - q^{6k+1} +q^{14k+3} 
    - q^{20k+5} - q^{42k+16}+q^{48k+20}  - q^{56k+26}+q^{62k+31} ~ \right] .  \nonumber           
   \end{eqnarray}  
($J_1(q)$ corresponds to sequence A108483.)  The coefficients in eq. (14) for $a_7 = 0,1,2,3,4$ and 6 are then            
            \begin{eqnarray*}
             {\bf Z}^{(0)}  = \left( \begin{array}{r} 1 \\ 7 \\ 35 \\ 12 \\  12 \\ -7 \\36 \\ -167 \\   \vdots~ \end{array} \right);~
             {\bf Z}^{(1)}  = \left( \begin{array}{r} 1  \\  14  \\  20 \\ 34  \\ -1 \\ 21 \\ -111 \\ 34 \\ \vdots ~\end{array} \right);~
             {\bf Z}^{(2)} = \left( \begin{array}{r}  2 \\ 14 \\ 31 \\ 7 \\ 44 \\-67 \\ 21 \\ -103 \\ \vdots ~\end{array} \right);~
             {\bf Z}^{(3)}  = \left( \begin{array}{r} 3 \\ 18  \\ 21 \\39 \\ -28 \\ 31  \\-80 \\ -73 \\  \vdots ~\end{array} \right);~
             {\bf Z}^{(4)}  = \left( \begin{array}{r} 5 \\ 16 \\ 37 \\ -2 \\  35 \\ -47 \\ -28\\ -117 \\  \vdots~ \end{array} \right);~
             {\bf Z}^{(6)} = \left( \begin{array}{r}  11 \\ 13 \\ 39 \\ 14 \\ 0 \\ -63 \\ -1 \\ -164 \\ \vdots ~\end{array} \right).               
            \end{eqnarray*}
 
 \section {Conclusion}   
  The reduction of the $(n +1)$-dimensional matrix in eq. (1) to a smaller, $(k+1)$-dimensional one for $p(7k+a_7)$, 
  as well as for $p(5k+a_5)$ and $p(25k+a_{25})$, used results from Ramanujan's study of the 
  congruences of the partition function in ref. 3.  However, these reductions did not depend upon the existence of a 
  congruence, but rather they used the property that the product 
         \begin{eqnarray*}
        (q)_{\infty} ( e^{2 \pi i/N} q)_{\infty} ( e^{4 \pi i/N} q)_{\infty}\cdots (e^{(N-1)2 \pi i/N}q)_{\infty}
        \end{eqnarray*} 
  is a power-series expansion in $q^N$.    Thus, for any $N$, with $\omega = e^{2 \pi i/N}$, we can write, 
       \begin{eqnarray}
     \sum_{n=0}^{\infty} p(n) q^n = \frac{1}{(q)_{\infty} } \times 
     \frac{ (\omega q)_{\infty}\cdots  (\omega^{N-1}q)_{\infty} }{ (\omega  q)_{\infty} \cdots  (\omega^{N-1}q)_{\infty} } 
     = \left( \sum_{p=0}^{N-1} \sum_{k=0}^{\infty} Z_k^{(p)} q^{kN+p} \right) \left(  \sum_{k=0}^{\infty} D_k q^{kN} \right)^{-1},
     \end{eqnarray}   
with $D_0 =1$ and $Z_0^{(a)}=p(a)$ .  We have then     
      \begin{eqnarray}
          p(kN+a) &=&  \left| \begin{array}{ccccc}
             1 & ~ & ~&~& p(a)\\
            D_1& 1& ~ & ~&Z_1^{(a)}\\
            D_2& D_1& 1 & ~&Z_2^{(a)}\\
            D_3& D_2&  D_1& ~& Z_3^{(a)}\\
             ~\vdots ~&~& ~& ~\ddots &  \vdots ~
            \end{array} \right|_{(k+1) \times (k+1)} .
           \end{eqnarray}
This is however only of practical use in calculating partition functions if compact expressions can be found whose expansions give the $D$ and $Z$ coefficients.


\begin{thebibliography}{99}
 \bibitem{Mal}  Malenfant, ``Finite, Closed-form Expressions for the Partition Function and for Euler, Bernoulli, and Stirling Numbers''. http://arxiv.org/abs/1103.1585
  \bibitem {integer} The On-Line Encyclopedia of Integer Sequences. http://oeis.org.
    \bibitem{Ram} Berndt and Ono, "Ramanujan's Unpublished Manuscript on the Partition and Tau Functions with
 Proofs and Commentary". http://www.math.wisc.edu/~ono/reprints/044.pdf.
\end{thebibliography}
\end{document}